\documentclass[11pt]{article}
\usepackage{amsmath}
\usepackage{amsthm}
\usepackage{bm}
\usepackage{multirow}

\usepackage{lipsum}
\usepackage{amsfonts}
\usepackage{graphicx}
\usepackage{epstopdf}
\usepackage{algorithm}
\usepackage{algorithmic}
\ifpdf
  \DeclareGraphicsExtensions{.eps,.pdf,.png,.jpg}
\else
  \DeclareGraphicsExtensions{.eps}
\fi

\newcommand{\TheTitle}{Multigrid waveform relaxation \\ for the time-fractional heat equation}

\title{{\TheTitle}
\footnote{Francisco J. Gaspar has received funding from the European Union's Horizon 2020 research and innovation programme under the Marie Sklodowska-Curie grant agreement No 705402, POROSOS. The work of Carmen Rodrigo is supported in part by the Spanish project FEDER /MCYT MTM2016-75139-R and the Diputaci\'on General de Arag\'on (Grupo consolidado PDIE).}}

  \author{
  Francisco J. Gaspar\footnote{CWI, Centrum Wiskunde and Informatica, Amsterdam, The Netherlands, (gaspar@cwi.nl)}
\and
  Carmen Rodrigo\footnote{IUMA and Applied Mathematics Department, University of Zaragoza, Spain, (carmenr@unizar.es)}
  }

\usepackage{amsopn}

\begin{document}

\maketitle

\begin{abstract}
In this work, we propose an efficient and robust multigrid method for solving the time-fractional heat equation. Due to the nonlocal property of fractional differential operators, numerical methods usually generate systems of equations for which
the coefficient matrix is dense. Therefore, the design of efficient solvers for the numerical simulation of these problems is a difficult task. We develop a parallel-in-time multigrid algorithm based on the waveform relaxation approach, whose application to time-fractional problems seems very natural due to the fact that the fractional derivative at each spatial point depends on the values of the function at this point at all earlier times. Exploiting the Toeplitz-like structure of the coefficient matrix, the proposed multigrid waveform relaxation method has a computational cost of $O(N M \log(M))$ operations, where $M$ is the number of time steps and $N$ is the number of spatial grid points. A semi-algebraic mode analysis is also developed to theoretically confirm the good results obtained. Several numerical experiments, including examples with non-smooth solutions and a nonlinear problem with applications in porous media, are presented.
\end{abstract}



\section{Introduction}\label{sec:1}

\setcounter{section}{1}
Fractional calculus has become increasingly popular in recent years due to their frequent appearance in various applications in fluid mechanics, signal processing, viscoelasticity, porous media flow, quantum mechanics, biology, medicine, physics and engineering, see~\cite{cushman, hall, hilfer2, lazarov,li_xu,metzler_klafter, purohit} for example.
In particular it has attracted much attention  within the natural and social sciences, since it can properly model phenomena dominated by memory effects~\cite{11hesthaven,10hesthaven} and problems exhibiting non-Markovian behavior in time.

A lot of effort has been focused in attempting to find robust and stable numerical 
and analytical 
methods for solving ordinary and partial differential equations of fractional order.
A wide growth in the number of numerical analysis papers studying differential equations with fractional-order derivatives have arisen in the last decade~\cite{BLY, lazarov, lazarov2, li_xu, lin_xu, mustapha, martin1}. Due to the nonlocal property of the fractional differential operator,  numerical methods usually generate systems of equations for which the coefficient matrix is dense. This is the main reason why most of  these problems have been traditionally solved by Gaussian elimination,  which requires a very high computational cost of $O(n^3)$ in addition to a high storage cost of $O(n^2)$, where $n$ is the total number of grid-points. Some efforts have been done to reduce this computational cost by approximating the coefficient matrix by a banded matrix \cite{GE_banded}, for example, obtaining a computational complexity of $O(n\log^2(n))$.
This is quite different from the integer differential operators, which typically yield sparse coefficient matrices that can be efficiently solved by fast iterative methods with 
$O(n)$ complexity. 
Therefore, the design of efficient solvers that reduce the computational cost is one of the challenges  for the numerical simulation of fractional PDEs.  For space-fractional PDEs, some efficient solvers, such as preconditioned Krylov subspace methods \cite{Pan_SISC,Wang_SISC} and multigrid methods \cite{Pang2012}, have already been proposed. The key is to take advantage of the Toeplitz-like structure of the coefficient matrix which arises from the discretization method proposed by Meerschaet and Tadjeran \cite{Meerschaert2006}. In this way, the storage requirements can be reduced to $O(N)$, and the matrix-vector multiplication can be done in $O(N \log(N))$ operations by using the fast Fourier transform (FFT). Recently, a fast solver based on a geometric multigrid method for nonuniform grids has been proposed  in~\cite{fractional_xiaozhe}. The approach is based on the use of H-Matrices to approximate the dense matrices. Regarding time-fractional PDEs, the coefficient matrix usually has an $M \times M$  block lower triangular Toeplitz structure, with each block of size $N \times N$, where $N$ is the number of spatial grid-points and $M$ the number of time levels. 
A fast direct method taking advantage of the Toeplitz structure of the coefficient matrix is proposed in \cite{GE_time_fractional} with a complexity of $O(NM\log^2(M))$.
As an alternative, in \cite{Zhang2014} the authors proposed the use of alternating direction implicit schemes (ADI) with a computational complexity of $O(NM^2)$. An approximate inversion method with $O(N M \log(M))$ has been recently proposed in \cite{Lin2016}.
Their idea is to approximate the coefficient matrix by a block $\varepsilon$-circulant matrix,
which can be block diagonalized by FFT. To solve the resulting complex block system, the authors use a multigrid method. Our main contribution is to propose an efficient and robust multigrid method based on the waveform relaxation approach to solve the time-fractional heat equation. 
Exploiting the Toeplitz-like structure of the coefficient matrix, the computational complexity of the proposed method is $O(N M \log(M))$ with a storage
requirement of $O(NM)$, being only $O(M)$ for the storage of the coefficient matrix. 
Opposite to the method introduced in \cite{Lin2016}, the algorithm proposed here is directly applied to the original discretization of the problem, and also is better suited for nonlinear problems.
We wish to emphasize that the proposed method 
is parallel-in-time in contrast to the classical sequential time-integration methods based on time-stepping.

Waveform relaxation methods consist of continuous-in-time iterative algorithms for solving large systems of ordinary differential equations (ODEs).
Their application to the solution of parabolic partial differential equations is based on the numerical method of lines, in which the spatial derivatives are replaced by discrete analogues obtaining a large system of ODEs, which is solved by standard iterative methods.
The requirement of extra storage for unknowns represents a classical disadvantage of waveform methods, however in our case this is not a drawback anymore since the time-fractional PDEs also need the solutions in previous time-steps to be stored. Since the waveform relaxation method is based on the numerical method of lines, it is not clear how to combine it with techniques such as dynamic grid adaptation, although recently some efforts have been carried out to combine parallel-in-time techniques with moving meshes~\cite{Falgout2016, Haynes_SISC}.
The convergence of the waveform relaxation methods was studied by Miekkala and Nevanlinna~\cite{miekkala}, who showed that the convergence could be too slow for the waveform relaxation to be competitive with standard time-stepping methods. Recently, some authors have investigated the convergence of some waveform relaxation methods for solving fractional differential equations~\cite{jiang_ding}. We wish to point out that for time-fractional PDEs, the fact of that each spatial grid point at a fixed time is connected to all the values of the previous time steps makes the application of waveform relaxation methods to these problems very natural.

Multigrid methods (see~\cite{Stu_Tro, TOS01, Wess} for an introduction) are often used for the convergence acceleration of iterative methods, although they have a wider use and significance than just being acceleration techniques. These methods are among the most efficient methods for solving large algebraic systems arising from discretizations of partial differential equations, with optimal computational complexity, due to their ability to handle different scales present in the problem.  Here, we propose the application of a multigrid approach based on the waveform relaxation method for solving time-fractional differential equations. This method combines the very fast multigrid convergence with the high parallel efficiency of waveform relaxation. Basically, it consists of applying a red-black zebra-in-time line relaxation together with a coarse-grid correction procedure based on coarsening only in the spatial dimension. Note that there is no coarsening in time in such a multigrid method, and the time is kept continuous. In this way, the coarsest grid is composed of only one spatial grid-point and all the corresponding points in time. The multigrid waveform relaxation was firstly developed by Lubich and Ostermann in~\cite{lubich_ostermann}, who showed that the basic waveform relaxation process can be accelerated by using the multigrid idea. Their work is based on the application of multigrid (in space) directly to the evolution equation. Since its introduction, this approach has been successfully applied to a variety of parabolic problems~\cite{horton_vandewalle, oosterlee_wesseling, computing, vandewalle_piessens, vandewalle_piessens3, vandewalle_piessens2}, but never within a fractional context.

Local Fourier analysis (LFA) or local mode analysis~\cite{Bra77, Bra94, TOS01, Wess, Wie01} has become a very useful predictive tool for the analysis of the convergence of multigrid methods. The idea of this analysis is to focus on the local character of the operators involved in the multigrid algorithm, and to analyze their behavior on a basis of complex exponential functions. However, the failure of this analysis for the prediction of the multigrid convergence for convection-dominated or parabolic problems has been observed by different authors~\cite{Brandt_singular, sama, OGWW}.  To overcome this difficulty, a semi-algebraic mode analysis (SAMA) was proposed in~\cite{sama}. This analysis, which is essentially a generalization of the classical local mode analysis, combines standard LFA with algebraic computation that accounts for the non-local character of the operators. It is clear that this is the approach that we should consider for the analysis of the multigrid waveform relaxation method for the time-fractional diffusion problem dealt with in this work. Notice the non-local character of this differential operator in time. Finally, we wish to emphasize that the proposed multigrid waveform relaxation method, as well as the semi-algebraic mode analysis for the study of its convergence, give rise an efficient solution strategy for the time-fractional heat equation, which seems a very natural way to deal with this problem.

The remainder of this paper is organized as follows. Section~\ref{sec:2} is devoted to introduce the considered one-dimensional model problem and its discretization. The proposed multigrid waveform relaxation method for its solution is described in Section~\ref{sec:3}. Next, the semi-algebraic mode analysis used for studying the convergence of this algorithm is explained in Section~\ref{sec:4}, together with some analysis results showing its suitability for the prediction of the behavior of the multigrid method. In Section~\ref{sec:5} the computational complexity of the proposed algorithm is discussed. After that, Section~\ref{sec:6} focuses on the generalization of the proposed methodology for a two-dimensional model problem. Finally, in Section~\ref{sec:7}, we illustrate the good behavior of the multigrid waveform relaxation method for solving the time-fractional diffusion problems considered in this work, by means of three different test problems, which include a nonlinear model problem with applications in porous media. Conclusions are drawn in Section~\ref{sec:8}.

\section{Model problem and discretization}\label{sec:2}
\setcounter{section}{2}

We consider the time-fractional heat equation, arising by replacing the first-order time derivative with the Caputo derivative of order $\delta$, where $0<\delta<1$. In this section, we restrict ourselves to the one-dimensional case for simplicity in the presentation. Therefore, we can formulate our model problem as the following initial-boundary value problem,
\begin{eqnarray}
D_t^{\delta} u - \frac{\partial^2 u}{\partial x^2} = f(x,t),\quad 0<x<L,\; t>0,\label{model_IVP_1}\\
u(0,t) = 0,\; u(L,t)=0,\quad t>0,\label{model_IVP_2}\\
u(x,0) = g(x), \quad 0\leq x\leq L.\label{model_IVP_3}
\end{eqnarray}
As mentioned above, $D_t^{\delta}$ denotes the Caputo fractional derivative, defined as follows~\cite{Diethelm,martin1}
\begin{equation}\label{caputo}
D_t^{\delta} u (x,t):= \left[J^{1-\delta}\left(\frac{\partial u}{\partial t}\right)\right](x,t),\quad 0\leq x\leq L,\; t>0,
\end{equation}
where $J^{1-\delta}$ represents the Riemann-Liouville fractional integral operator of order $1-\delta$, given by
\begin{equation}\label{RiemanLiouville}
\left(J^{1-\delta}u\right)(x,t):=\left[\frac{1}{\Gamma(1-\delta)}\int_{0}^{t} (t-s)^{-\delta}u(x,s)ds\right],\quad 0\leq x\leq L,\; t>0,
\end{equation}
where $\Gamma$ is the Gamma function~\cite{gamma}. 

Model problem~\eqref{model_IVP_1}-\eqref{model_IVP_3} is discretized on a uniform rectangular mesh $G_{h,\tau} = G_h \times G_{\tau},$ with
\begin{eqnarray}
G_{h} &=& \left\{x_n = nh,\, n=0,1,\ldots,N+1\right\},\label{spatial_mesh}\\
G_{\tau} &=& \left\{t_m= m\tau,\, m=0,1,\ldots,M\right\},\label{temporal_mesh}
\end{eqnarray}
where $h=\displaystyle \frac{L}{N+1}$, $\tau = \displaystyle\frac{T}{M}$ with $T$ the final time and $N+1$ and $M$ positive integers representing the number of subdivisions in the spatial and temporal intervals, respectively. We denote as $u_{n,m}$ the nodal approximation to the solution at each grid point $(x_n,t_m)$.\\
In the nodal points, the Caputo fractional derivative $D_t^{\delta}u$ can be written as follows 
\begin{equation}\label{caputo_nodal}
D_t^{\delta}\, u(x_n,t_m) = \frac{1}{\Gamma(1-\delta)}\sum_{k=0}^{m-1}\int_{t_k}^{t_{k+1}}(t_m-s)^{-\delta}\,\frac{\partial u(x_n,s)}{\partial t}ds,
\end{equation}
and it is approximated by using the well-known L1 scheme~\cite{Oldham_Spanier} which uses $\displaystyle \frac{\partial u(x_n,s)}{\partial t}\approx\frac{u_{n,k+1}-u_{n,k}}{\tau},\; t_k\leq s\leq t_{k+1}$ to obtain
\begin{eqnarray}
D_M^{\delta}u_{n,m}&:=&\frac{1}{\Gamma(1-\delta)}\sum_{k=0}^{m-1}\frac{u_{n,k+1}-u_{n,k}}{\tau} \int_{t_k}^{t_{k+1}}(t_m-s)^{-\delta}ds \label{L1_scheme1}
\\&=&\frac{\tau^{-\delta}}{\Gamma(2-\delta)}\left[d_1 u_{n,m} - d_m u_{n,0} + \sum_{k=1}^{m-1}(d_{k+1}-d_k)u_{n,m-k}\right],\label{L1_scheme2}
\end{eqnarray}
by defining $d_k = k^{1-\delta}-(k-1)^{1-\delta},\; k\geq 1$.\\
Regarding the diffusive term, we use standard central finite differences to approximate the spatial derivatives. Summarizing, we treat with the following discrete problem
\begin{eqnarray}
D_M^{\delta} u_{n,m} - \frac{u_{n+1,m}-2u_{n,m}+u_{n-1,m}}{h^2} = f(x_n,t_m), \; &1\leq n\leq N, \; 1\leq m\leq M,& \label{discrete_model_IVP_1} \\
u_{0,m} = 0,\; u_{N+1,m}=0, &0<m\leq M,& \label{discrete_model_IVP_2} \\
u_{n,0} = g(x_n), &0\leq n\leq N+1.& \label{discrete_model_IVP_3}
\end{eqnarray}

\section{Multigrid waveform relaxation in 1D}\label{sec:3}

\setcounter{section}{3}

For solving time dependent partial differential equations, the multigrid waveform relaxation method uses the numerical method of lines, replacing any spatial derivative by discrete formulas (obtained by the finite element, finite difference or finite volume methods) in the discrete spatial domain. Thus, the PDE is transformed to a large set of ordinary differential equations. In our case, that is, considering time-fractional derivatives of order $\delta$, we have
\begin{equation}\label{system_ODEs}
D_t^{\delta} u_h(t) + A_hu_h(t) = f_h(t),\; u_h(0) = g_h,\; t>0,
\end{equation}
where $u_h$ and $f_h$ are functions of time $t$ defined on a discrete spatial mesh, and $A_h$ is the discrete approximation in space of the continuous operator defining the considered PDE. Since discrete operators are usually represented by matrices and grid-functions by vectors, in the following we will use either the terminology of discrete differential operators and grid-functions or that of matrices and vectors.
Next step is the solution of the large system of ODEs by an iterative algorithm. For instance, if we consider a splitting of the spatial discrete operator $A_h = M_h-N_h$, one step of the iterative scheme for~\eqref{system_ODEs} can be written as
\begin{equation}\label{GS_ODEs}
D_t^{\delta} u_h^{k}(t) + M_hu_h^k(t) = N_h u_h^{k-1}(t) + f_h(t),\; u_h^k(0) = g_h, \; \hbox{for} \; k\geq 1,
\end{equation}
where $u_h^k(t)$ denotes the approximation obtained at iteration $k$.
The initial iterate $u_h^0(t)$ is defined along the whole time-interval, being natural to choose a constant initial approximation equal to the initial condition in~\eqref{system_ODEs}, that is, $u_h^0(t) = g_h,\; t>0$. \\
In this work, for the one-dimensional problem, we will consider a red-black Gauss-Seidel iteration which consists of a two-stage procedure, given by
\begin{eqnarray}\label{red-black}
&& D_t^{\delta}u_{n}^k(t) +\frac{2}{h^2}u_{n}^k(t) = \frac{1}{h^2}\left(u_{n-1}^{k-1}(t)+u_{n+1}^{k-1}(t) \right)+f_n(t),\; \hbox{if} \; n\; \hbox{is even},\\
&& D_t^{\delta}u_{n}^k(t) -\frac{1}{h^2}\left(u_{n-1}^{k}(t)-2u_{n}^k(t)+u_{n+1}^{k}(t) \right) = f_n(t),\; \hbox{if} \; n\; \hbox{is odd},
\end{eqnarray}
that is, first the even points in space are visited and after that we solve the unknowns in the grid points with odd numbering.

To accelerate the convergence of the red-black Gauss-Seidel waveform relaxation, a coarse-grid correction process based on a coarsening procedure only in the spatial dimension is performed, resulting the so-called linear multigrid waveform relaxation algorithm~\cite{Vandewalle_book}. This method consists essentially in the standard multigrid algorithm but applied to systems of ODEs as the one in~\eqref{system_ODEs}. Considering the standard full-weighting restriction and the linear interpolation as transfer-grid operators, the algorithm of the multigrid waveform relaxation (WRMG) is given in Algorithm~\ref{wrmg_1}.
\begin{algorithm}[htb]
\caption{ \textbf{: Multigrid waveform relaxation: ${\mathbf{u_{h}^k(t) \rightarrow u_{h}^{k+1}(t)}}$}}\label{wrmg_1}
\vspace{0.5cm}
\begin{algorithmic}
\IF{we are on the coarsest grid-level (with spatial grid-size given by $h_0$)}
\STATE
\STATE $D_t^{\delta} u_{h_0}^{k+1}(t) + A_{h_0}u_{h_0}^{k+1}(t)= f_{h_0}(t)$ \hspace{0.9cm} {Solve with a direct or fast solver.}
\STATE
\ELSE
\STATE {
\begin{tabular}{lr}
\\ [-1.5ex]
$\overline{u}_{h}^k(t)=S_{h}^{\nu_1}(u_{h}^k(t))$ & \hspace{-3.5cm}\textbf{(Pre-smoothing)}\\ & \hspace{-3.5cm}$\nu_1$ steps of the \textbf{red-black waveform relaxation}.\\
\\ [-1.5ex]
$\overline{r}_{h}^k(t)=f_{h}(t) - (D_t^{\delta} + A_{h})\, \overline{u}_{h}^k(t)$ & \hspace{-3.5cm}Compute the defect. \\
\\ [-1.5ex]
$\overline{r}_{2h}^k(t)=I_h^{2h}\, \overline{r}_{h}^k(t)$ & \hspace{-3.5cm}Restrict the defect. \\
\\ [-1.5ex]
$(D_t^{\delta}+A_{2h}) \widehat{e}_{2h}^k(t) = \bar{r}_{2h}^k(t),\; \widehat{e}_{2h}^k(0) = 0$ & \hspace{-3.5cm}Solve the defect equation \\ [1.2ex]
& \hspace{-3.5cm} on $G_{2h}$ by performing $\gamma \ge 1$ cycles of WRMG. \\
\\ [-1.5ex]
$\widehat{e}_{h}^k(t) = I_{2h}^h \, \widehat{e}_{2h}^k(t) $ & \hspace{-3.5cm}Interpolate the correction. \\
\\ [-1.5ex]
$\overline{u}_{h}^{k+1}(t) = \overline{u}_{h}^k(t) + \widehat{e}_{h}^k(t)$ & \hspace{-3.5cm}Compute a new approximation. \\
\\ [-1.5ex]
$u_{h}^{k+1}(t)=S_{h}^{\nu_2}(\overline{u}_{h}^{k+1}(t))$ & \hspace{-3.5cm}\textbf{(Post-smoothing)}\\ & \hspace{-3.5cm}$\nu_2$ steps of the \textbf{red-black waveform relaxation}. \\
\\ [-1.5ex]
\end{tabular}}
\ENDIF
\end{algorithmic}
\vspace{0.5cm}
\end{algorithm}

After discretizing in time, that is, replacing the differential operator $D_t^{\delta}$ by $D_M^{\delta}$, the previous algorithm can be interpreted as a space-time multigrid method with coarsening only in space.
Thus, the red-black Gauss-Seidel waveform relaxation can be seen as a zebra-in-time line relaxation, and standard full-weighting restriction and linear interpolation in space are considered for the data transfer between the levels in the multigrid hierarchy.
Thus, the whole multigrid waveform relaxation combines a zebra-in-time line relaxation with a standard semi-coarsening strategy only in the spatial dimension.
\section{Semi-algebraic mode analysis in 1D}\label{sec:4}

\setcounter{section}{4}

The analysis that we perform here is based on an exponential Fourier mode analysis or local Fourier analysis technique only in space and an exact analytical approach in time. This kind of semi-algebraic mode analysis was introduced for the first time in~\cite{sama}, where the authors mainly study the convergence of multigrid methods on space-time grids for parabolic problems. Furthermore, they extend the application of this analysis to non-parabolic problems like elliptic diffusion in layered media and convection diffusion. The main idea of this analysis is to study the evolution of the spatial Fourier modes over time. This semi-algebraic analysis provides very accurate predictions of the performance of multigrid methods, and indeed, it can be made rigorous if appropriate boundary conditions are considered.
Next, we describe the basics of this analysis. Although in~\cite{sama} the authors give a detailed explanation, here we present a slightly different description of this analysis.

\subsection{Basics of the analysis}\label{sec:4_1}

It is well-known that LFA assumes the formal extension to all multigrid components to an infinite grid, neglecting the boundary conditions, and considers discrete linear operators with constant coefficients. Therefore, we define the following infinite grid:
\begin{equation}\label{infinite_grid}
{\mathcal G}_h = \left\{x_n = n\,h,\; n\in {\mathbb Z}\right\},
\end{equation}
where $h$ is the spatial discretization step. For a fixed $t$, any discrete grid-function $u_h(\cdot,t)$ defined on ${\mathcal G}_h$ can be written as a formal linear combination of the so-called Fourier modes given by $\varphi_h(\theta,x) = e^{\imath\theta x}$, where $\theta\in\Theta_h=(-\pi/h,\pi/h]$, that is,
\begin{equation}\label{linear_combination}
u_h(x,t) = \sum_{\theta\in\Theta_h} c_{\theta}(t) \varphi_h(\theta,x), \quad x\in {\mathcal G}_h.
\end{equation}
Notice that coefficients $c_{\theta}(t)$ depend on the time variable. The Fourier modes, which generate the so-called Fourier space ${\mathcal F}({\mathcal G}_h) = \left\{\varphi_h(\theta,\cdot),\; \theta\in \Theta_h\right\}$, result to be eigenfunctions of any discrete operator with constant coefficients defined formally on ${\mathcal G}_h$. For instance, for the discrete operator $A_h = \displaystyle\frac{1}{h^2}\left[-1 \; 2\; -1\right]$, considered in discrete model problem~\eqref{discrete_model_IVP_1}, it is fulfilled that
$$A_h \varphi_h (\theta,\cdot) = \widehat{A}_h(\theta)\,\varphi_h(\theta,\cdot),$$
where $$\widehat{A}_h(\theta) = \displaystyle\frac{4}{h^2}\sin\left(\frac{\theta h}{2}\right)$$
is the Fourier representation of operator $A_h$, which is also called the Fourier symbol of $A_h$.

The aim of the local Fourier analysis is to analyze how the operators involved in the multigrid algorithm act on such Fourier modes. We can study how efficiently the relaxation process eliminates the high-frequency components of the error, through a smoothing analysis, or how the two-grid operator acts on the Fourier space, through a two-grid analysis. \\
First, we proceed to explain the smoothing analysis for a standard relaxation procedure. After that, we describe the analysis for the coarse-grid correction operator, and finally we combine both analysis in order to perform a complete study of the two-grid cycle. For this purpose, we need to distinguish high- and low-frequency components. This classification is done with respect to the coarsening strategy, which is chosen as standard coarsening, that is, the step size is double on the coarse grid, which is denoted by ${\mathcal G}_{2h}$. Remind that in a typical multigrid waveform relaxation procedure the coarsening applies only in the spatial domain. The space of low frequencies is defined as $ \Theta_{2h} = (-\pi/2h,\pi/2h]$, and the high-frequencies are given by $ \Theta_h\backslash \Theta_{2h}$.\\

\noindent {\bf Smoothing analysis.} We describe the semi-algebraic smoothing analysis for a standard relaxation procedure based on a decomposition of the spatial discrete operator $A_h$ as $A_h = M_h - N_h$. Denoting $e_h^k(\cdot,t)$ and $e_h^{k-1}(\cdot,t)$ the error grid-functions at the $k$ and $k-1$ iterations of this procedure, an iteration of this waveform relaxation method is given by
\begin{equation}\label{waveform_error}
D_t^{\delta}e_h^{k}(x,t) + M_h e_h^k(x,t) = N_h e_h^{k-1}(x,t), \; \hbox{for} \; k\geq 1, \;\hbox{and} \; x\in{\mathcal G}_h, \; t>0,
\end{equation}
with initial condition $e_h^k(x,0) = 0,\; x\in{\mathcal G}_h$.\\
From~\eqref{linear_combination}, we can write the error at $j$ iteration, $e_h^j(x,t)$, in the following way,
\begin{equation}\label{error_combinations_j}
e_h^j(x,t) = \sum_{\theta\in\Theta_h} c_{\theta}^j(t) \varphi_h(\theta,x), \quad x\in {\mathcal G}_h,\; t>0,\\
\end{equation}
and then by using that $\varphi_h(\theta,x)$ are eigenfunctions of operators $M_h$ and $N_h$ (that is, $M_h\varphi_h(\theta,x) = \widehat{M}_h(\theta) \varphi_h(\theta,x)$ for example, where $\widehat{M}_h(\theta)$ is the symbol of $M_h$), it follows for each frequency $\theta\in\Theta_h$ that
\begin{equation}\label{GS_error_theta}
D_t^{\delta}c_{\theta}^k(t) + \widehat{M}_h(\theta)c_{\theta}^k(t) = \widehat{N}_h(\theta) c_{\theta}^{k-1}(t),\; \hbox{for} \; k\geq 1,\; t>0.
\end{equation}
Considering the discretization of $D_t^{\delta}$ on the uniform grid $G_{\tau}$, $D_M^{\delta}$, defined in~\eqref{L1_scheme2}, and denoting $\left( c_{\theta}^{k,1},\ldots, c_{\theta}^{k,M}\right)$ the approximation of $c_{\theta}^{k}(t)$ on grid $G_{\tau}$, we obtain the following relation
\begin{equation}\label{algebraic_system}
\left(\begin{array}{c} c_{\theta}^{k,1} \\  c_{\theta}^{k,2}\\ \vdots \\ c_{\theta}^{k,M} \end{array}\right) = \widetilde{{\mathcal M}}_{h,\tau}^{-1}(\theta) \widetilde{{\mathcal N}}_{h,\tau}(\theta)
\left(\begin{array}{c} c_{\theta}^{k-1,1} \\  c_{\theta}^{k-1,2}\\ \vdots \\ c_{\theta}^{k-1,M} \end{array}\right),
\end{equation}
where $\widetilde{{\mathcal N}}_{h,\tau} (\theta)= \hbox{diag}(\widehat{N}_h(\theta))$, and $$\widetilde{{\mathcal M}}_{h,\tau}(\theta) = \left(\begin{array}{cccc} r_1 +\widehat{M}_h(\theta) & 0 & \cdots & 0 \\  r_2 & r_1 + \widehat{M}_h(\theta) & \cdots & 0 \\ \vdots & \ddots & \ddots & \vdots \\ r_M & \cdots & r_2 & r_1 + \widehat{M}_h(\theta) \end{array}\right),$$ with $r_i = \displaystyle \frac{\tau^{-\delta}}{\Gamma(2-\delta)} (d_i-d_{i-1}), i = 1, \ldots, M$, assuming $d_0 = 0$.\\
Denoting $\widetilde{{\mathcal S}}_{h,\tau} (\theta)=\widetilde{{\mathcal M}}_{h,\tau}^{-1}(\theta) \widetilde{{\mathcal N}}_{h,\tau}(\theta)$, we can define the smoothing factor of the relaxation procedure as follows
\begin{equation}\label{smoothing_factor}
\mu = \sup_{\Theta_h\backslash\Theta_{2h}} \left(\rho\left(\widetilde{{\mathcal S}}_{h,\tau}(\theta)\right)\right).
\end{equation}

\noindent {\bf Coarse-grid correction analysis.} We now proceed to explain the analysis of the coarse-grid correction method. An error $e_h^k$ is transformed by this method as $e_h^{k+1} = C_h^{2h} e_h^k$, where  $C_h^{2h} = I_h - I_{2h}^h (D_t^{\delta} + A_{2h})^{-1} I_h^{2h} (D_t^{\delta} + A_h)$ is the coarse-grid correction operator. Here $I_h$ is the identity operator, $D_t^{\delta}+A_h$ and $D_t^{\delta} + A_{2h}$ are the fine- and coarse-grid operators, respectively, and $I_{2h}^h$, $I_h^{2h}$ are transfer operators from coarse to fine grids and vice versa. \\
As we have chosen standard coarsening, the fine-grid Fourier mode $\varphi_h(\theta,x)$ when injected
into the coarse grid, aliases with the coarse-grid Fourier mode $\varphi_{2h}(2\theta,x)$. Thus, for any low-frequency $\theta^0 \in \Theta_{2h}$, we define the high-frequency $\theta^1 = \theta^0 - {\rm sign}(\theta^0) \pi/h$. Taking this into account, the Fourier space is decomposed into two-dimensional subspaces, known as 2h-harmonics (see~\cite{TOS01,Wie01} for more details):
\[
{\mathcal F}^2(\theta) = {\rm span} \{ \varphi_h(\theta^0,\cdot),  \varphi_h(\theta^1,\cdot) \}, \quad \theta = \theta^0 \in
\Theta_{2h}.
\]
The coarse-grid correction operator $C_h^{2h}$ leaves the two-dimensional subspace of harmonics
${\mathcal F}^2(\theta^0)$ invariant for an arbitrary Fourier frequency $\theta^0 \in \Theta_{2h}$.
Let us define for any $\theta^0 \in \Theta_{2h}$ the vector
$\boldsymbol{\varphi}_h(\theta^0,\cdot) = (\varphi_h(\theta^0,\cdot),\varphi_h(\theta^1,\cdot))$. As the error at the iteration k can be written as $e_h^k(x,t) = \sum_{\theta \in \Theta_{2h}} {\mathbf c}_{\theta}^k(t) \boldsymbol{\varphi}_h(\theta,x)^T$, with ${\mathbf c}_{\theta}^k(t) = (c_{\theta^0}^k(t),c_{\theta^1}^k(t))$, the error at the iteration $k+1$ after application of the coarse-grid correction method is given by $\sum_{\theta \in \Theta_{2h}}  \widehat{C}_h^{2h}(\theta) {\mathbf c}_{\theta}^k(t) \boldsymbol{\varphi}_h(\theta,\cdot)^T$, where $\widehat{C}_h^{2h}(\theta)$ is a $2\times2$ matrix given by the expression
\[
\widehat{C}_h^{2h}(\theta) = I_2- \widehat{I}_{2h}^h(\theta) (D_t^{\delta} + \widehat{A}_{2h}(\theta))^{-1} \widehat{I}_h^{2h}(\theta) (D_t^{\delta} + \widehat{A}_h(\theta)),
\]
where $I_2$ is the $2\times 2$ identity matrix, and $\widehat{A}_h(\theta), \widehat{A}_{2h}(\theta), \widehat{I}_{2h}^h(\theta), \widehat{I}_h^{2h}(\theta)$ denote the symbols of the fine- and coarse-grid spatial operators, the prolongation operator, and the restriction operator, respectively. The Fourier symbol of the fine-grid operator is given by $\widehat{A}_h(\theta) = {\rm diag}(\widehat{A}_h(\theta^0),\widehat{A}_h(\theta^1))$, and the symbol of the coarse grid operator by $\widehat{A}_{2h}(\theta)$. The Fourier symbols of the prolongation and restriction operators for $\theta = \theta^0\in\Theta_{2h}$ are given by
\[
\widehat{I}_{2h}^h(\theta) =
\left(
\begin{array}{c}
\widehat{I}_{2h}^h(\theta^0) \\
\widehat{I}_{2h}^h(\theta^1)
\end{array}
\right),
\qquad \qquad
\widehat{I}_{h}^{2h}(\theta) = (\widehat{I}_{h}^{2h}(\theta^0) ,\widehat{I}_{h}^{2h}(\theta^1)).
\]
Let us suppose that the error at the iteration $k$ is given by $${\mathbf c}_{\theta}^k(t) \boldsymbol{\varphi}_h(\theta,\cdot)^T = c_{\theta^0}^k(t) \varphi_h(\theta^0,\cdot) + c_{\theta^1}^k(t) \varphi_h(\theta^1,\cdot).$$ By considering the discretization of $D_t^{\delta}$ on the uniform grid $G_{\tau}$, $D_M^{\delta}$, defined in~\eqref{L1_scheme2}, we obtain that the error after application of the coarse-grid correction is given by $\widetilde{{\mathcal C}}_{h,\tau}^{2h}(\theta) {\mathbf c}_{\theta}^k(t) \boldsymbol{\varphi}_h(\theta,\cdot)^T$, with $\widetilde{{\mathcal C}}_{h,\tau}^{2h}(\theta)$ a $2M \times 2M$  matrix, given by
\[
\widetilde{{\mathcal C}}_{h,\tau}^{2h}(\theta) = I_{2M}- \widetilde{{\mathcal I}}_{2h}^h(\theta) (\widetilde{{\mathcal A}}_{2h,\tau}(\theta))^{-1} \widetilde{{\mathcal I}}_h^{2h}(\theta)  \widetilde{{\mathcal A}}_{h,\tau}(\theta).
\]
Here, $I_{2M}$ is the identity matrix of order $2M$, $\widetilde{{\mathcal A}}_{h,\tau}(\theta)$ is the $2M\times 2M$ matrix
\[
\widetilde{{\mathcal A}}_{h,\tau}(\theta)) = \left(
\begin{array}{cc}
\widetilde{{\mathcal A}}_{h,\tau}(\theta^0) & 0 \\
0 & \widetilde{{\mathcal A}}_{h,\tau}(\theta^1)
\end{array}
\right), \quad \theta = \theta^0\in \Theta_{2h},
\]
where for $\alpha=0,1$,
\begin{equation}\label{fine_grid_operator_symbol}
\widetilde{{\mathcal A}}_{h,\tau}(\theta^{\alpha}) = \left(\begin{array}{cccc} r_1 + \widehat{A}_h(\theta^{\alpha}) & 0 & \cdots & 0 \\  r_2 & r_1 + \widehat{A}_h(\theta^{\alpha}) & \cdots & 0 \\ \vdots & \ddots & \ddots & \vdots \\ r_M & \cdots & r_2 & r_1 + \widehat{A}_h(\theta^{\alpha}) \end{array}\right),
\end{equation}
with $\widehat{A}_h(\theta^{\alpha})$ the symbol of the fine-grid spatial operator, and $r_i = \displaystyle \frac{\tau^{-\delta}}{\Gamma(2-\delta)} (d_i-d_{i-1})$,  for $i=1,\ldots, M$, assuming $d_0 = 0$. \\
About the restriction and interpolation, $\widetilde{{\mathcal I}}_h^{2h}(\theta)$ is the matrix $M\times 2M$
\[
\widetilde{{\mathcal I}}_h^{2h}(\theta) = \left[  \widehat{I}_{h}^{2h}(\theta^0) I_M ,\widehat{I}_{h}^{2h}(\theta^1) I_M \right],
\]
and $\widetilde{{\mathcal I}}_{2h}^{h}(\theta)$ is the matrix $2M\times M$
\[
\widetilde{{\mathcal I}}_{2h}^{h}(\theta) = \left[  \widehat{I}_{2h}^{h}(\theta^0) I_M ,\widehat{I}_{2h}^{h}(\theta^1) I_M \right]^T.
\]
\noindent {\bf Two-grid analysis.} Combining the Fourier smoothing analysis and the Fourier coarse-grid correction analysis previously introduced, we perform the semi-algebraic two-grid analysis. The two-grid operator ${\mathcal T}_{h,\tau}^{2h}$ is defined as ${\mathcal T}_{h,\tau}^{2h} = {\mathcal S}_{h,\tau}^{\nu_2} {\mathcal C}_{h,\tau}^{2h} {\mathcal S}_{h,\tau}^{\nu_1}$, where ${\mathcal C}_{h,\tau}^{2h}$ is the coarse-grid operator, ${\mathcal S}_{h,\tau}$ a smoothing operator, and $\nu_1$, $\nu_2$ indicate the number of pre- and post-smoothing steps, respectively. \\

We remind that the  coarse grid correction operator ${\mathcal C}_{h,\tau}^{2h}$ leaves the two-dimen\-sional subspaces of harmonics ${\mathcal F}^2(\theta)$ invariant for an arbitrary Fourier frequency $\theta = \theta^0 \in \Theta_{2h}$. This same invariance property is true for the smoothers ${\mathcal S}_{h,\tau}$ considered in this work. Therefore, the two-grid operator ${\mathcal T}_{h,\tau}^{2h}$ also leaves the 2h-harmonic subspaces invariant. \\

Let us suppose that the error at the iteration $k$ is given by ${\mathbf c}_{\theta}^k(t) \boldsymbol{\varphi}_h(\theta,\cdot)^T = c_{\theta^0}^k(t) \varphi_h(\theta^0,\cdot) + c_{\theta^1}^k(t) \varphi_h(\theta^1,\cdot)$. By considering the discretization of $D_t^{\delta}$ on the uniform grid $G_{\tau}$, $D_M^{\delta}$, defined in~\eqref{L1_scheme2}, we obtain that the error after application of the two-grid method is given by $\widetilde{{\mathcal T}}_{h,\tau}^{2h}(\theta) {\mathbf c}_{\theta}^k(t) \boldsymbol{\varphi}_h(\theta,\cdot)^T$, with $\widetilde{{\mathcal T}}_{h,\tau}^{2h}(\theta)$ a $2M \times 2M$  matrix, given by
\[
\widetilde{{\mathcal T}}_{h,\tau}^{2h}(\theta) = \widetilde{{\mathcal S}}_{h,\tau}^{\nu_2}(\theta) (I_{2M}- \widetilde{{\mathcal I}}_{2h}^h(\theta) (\widetilde{{\mathcal A}}_{2h,\tau}(\theta))^{-1} \widetilde{{\mathcal I}}_h^{2h}(\theta)  \widetilde{{\mathcal A}}_{h,\tau}(\theta))\widetilde{{\mathcal S}}_{h,\tau}^{\nu_1}(\theta).
\]
If the chosen smoother is an iterative method which does not couple frequencies, then $\widetilde{{\mathcal S}}_{h,\tau}(\theta)$ is the $2M\times 2M$ matrix
\[
\widetilde{{\mathcal S}}_{h,\tau}(\theta) = \left(
\begin{array}{cc}
\widetilde{{\mathcal S}}_{h,\tau}(\theta^0) & 0 \\
0 & \widetilde{{\mathcal S}}_{h,\tau}(\theta^1)
\end{array}
\right),
\]
where for $\alpha=0,1$, $\widetilde{{\mathcal S}}_{h,\tau}(\theta^{\alpha})$ is given as previously. \\
In the case of a pattern waveform relaxation method, as the red-black waveform relaxation, it is well-known that the smoother couples frequencies but leaves invariant the two-dimensional subspaces ${\mathcal F}^2(\theta)$.
In particular, for the red-black waveform relaxation considered in this work, the symbol is given by $\widetilde{{\mathcal S}}_{h,\tau}(\theta) = \widetilde{{\mathcal S}}_{h,\tau}^{black}(\theta)\widetilde{{\mathcal S}}_{h,\tau}^{red}(\theta)$, where $\widetilde{{\mathcal S}}_{h,\tau}^{black}(\theta)$ and $\widetilde{{\mathcal S}}_{h,\tau}^{red}(\theta)$ are $2M\times 2M$ matrices coupling frequencies $\theta^0$ and $\theta^1$. More concretely,
$$\widetilde{{\mathcal S}}_{h,\tau}^{red}(\theta) = \frac{1}{2}\left(\begin{array}{cc} \widetilde{{\mathcal M}}_{h,\tau}^{-1}(\theta^0)\widetilde{{\mathcal N}}_{h,\tau}(\theta^0)+I_M & \widetilde{{\mathcal M}}_{h,\tau}^{-1}(\theta^1)\widetilde{{\mathcal N}}_{h,\tau}(\theta^1)-I_M \\ \widetilde{{\mathcal M}}_{h,\tau}^{-1}(\theta^0)\widetilde{{\mathcal N}}_{h,\tau}(\theta^0)-I_M & \widetilde{{\mathcal M}}_{h,\tau}^{-1}(\theta^1)\widetilde{{\mathcal N}}_{h,\tau}(\theta^1)+I_M \end{array}\right),$$
$$\widetilde{{\mathcal S}}_{h,\tau}^{black}(\theta) = \frac{1}{2}\left(\begin{array}{cc} \widetilde{{\mathcal M}}_{h,\tau}^{-1}(\theta^0)\widetilde{{\mathcal N}}_{h,\tau}(\theta^0)+I_M & -\widetilde{{\mathcal M}}_{h,\tau}^{-1}(\theta^1)\widetilde{{\mathcal N}}_{h,\tau}(\theta^1)+I_M \\ -\widetilde{{\mathcal M}}_{h,\tau}^{-1}(\theta^0)\widetilde{{\mathcal N}}_{h,\tau}(\theta^0)+I_M & \widetilde{{\mathcal M}}_{h,\tau}^{-1}(\theta^1)\widetilde{{\mathcal N}}_{h,\tau}(\theta^1)+I_M \end{array}\right),$$
where $I_M$ is the identity matrix of size $M\times M$, and $\widetilde{{\mathcal M}}_{h,\tau}(\theta)$ and $\widetilde{{\mathcal N}}_{h,\tau}(\theta)$ are given as explained in the smoothing analysis section, using that $M_h$ is the diagonal part of matrix $A_h$ as usual for a Jacobi-type relaxation. For a more detailed explanation of the semi-algebraic mode analysis for this smoother we refer to the reader to~\cite{sama}.\\

Finally, the convergence factor of the two-grid method, can be estimated as
\begin{equation}\label{two_grid_factor}
\rho = \sup_{\Theta_{2h}} \left(\rho\left(\widetilde{{\mathcal T}}_{h,\tau}^{2h}(\theta)\right)\right),
\end{equation}

\subsection{Analysis results}\label{sec:4_2}

\begin{figure}[htb]
\begin{center}
\includegraphics[width = 0.8\textwidth]{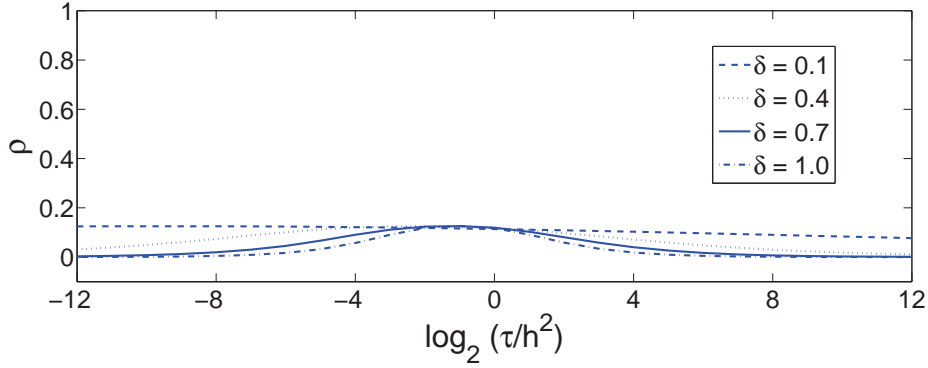}
\caption{Two-grid convergence factors predicted by the analysis for different values of parameter $\lambda = \tau/h^2$ from $2^{-12}$ to $2^{12}$ and different fractional orders $\delta$.}
\label{delta_lambda_star}
\end{center}
\end{figure}

\begin{figure}[htb]
\begin{center}
\includegraphics[width = 0.8\textwidth]{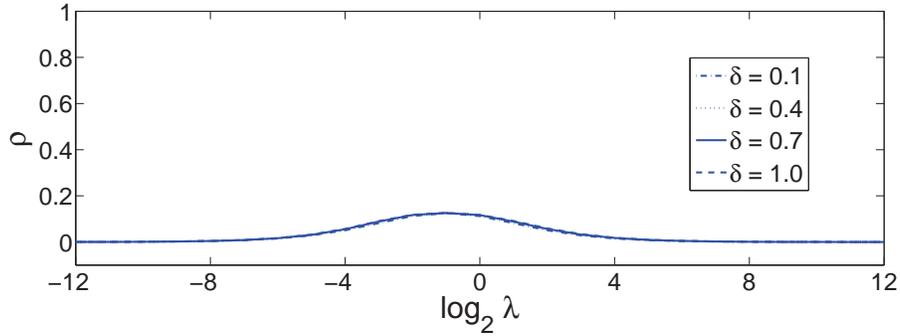}
\caption{Two-grid convergence factors predicted by the analysis for different values of parameter $\lambda = \tau^{\delta} \Gamma(2-\delta)/h^2$ and various fractional orders $\delta$.}
\label{delta_lambda}
\end{center}
\end{figure}

This section is focused on the analysis of the robustness of the proposed multigrid waveform relaxation method for the considered problem. When studying the multigrid convergence for the standard heat equation, it is well-known that parameter $\tau/h^2$ describes the anisotropy in the operator, resulting the relevant parameter for its analysis, see~\cite{paper_schrodinger} for example.
However, as it can be observed in Figure~\ref{delta_lambda_star}, this parameter is not the important one for the time-fractional heat equation. In Figure~\ref{delta_lambda_star}, we depict the two-grid convergence factors provided by the semi-algebraic mode analysis for a range of values of parameter $\tau/h^2$ from $2^{-12}$ to $2^{12}$, for different fractional orders $\delta$. Only one smoothing step is considered, and the zebra-in-time smoother is used as previously described. It is clearly seen that, although the convergence rates are bounded by $0.2$ for all cases, we do not obtain a $\delta-$independent convergence for a fixed value of $\tau/h^2$.
In this case, the relevant parameter is $\lambda = \tau^{\delta} \Gamma(2-\delta)/h^2$, as shown in Figure~\ref{delta_lambda}, where it is observed that the obtained multigrid convergence becomes robust for any value of $\delta$ with respect to parameter $\lambda$. In this figure, the number of time-steps is chosen as $M = 32$.
Notice that, for any fixed value of $\delta$, the multigrid convergence is satisfactory for any value of parameter $\lambda$, which is very important for the global behavior of the method since this parameter will vary from grid-level to grid-level within the multigrid algorithm. The corresponding MATLAB function used to carry out the SAMA results in this figure is available as supplementary material.

\begin{figure}[htb]
\begin{center}
\includegraphics[width = 0.8\textwidth]{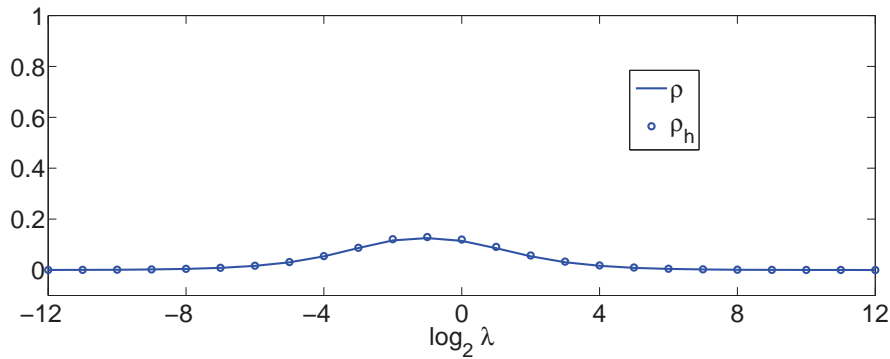}
\caption{Comparison between the two-grid convergence factors predicted by the analysis ($\rho$) and the asymptotic convergence factor of a $W-$cycle experimentally computed ($\rho_h$), for different values of parameter $\lambda = \tau^{\delta} \Gamma(2-\delta)/h^2$ and fractional order $\delta = 0.4$.}
\label{comparison}
\end{center}
\end{figure}

The results obtained by the semi-algebraic mode analysis match very accurately the real asymptotic convergence factors experimentally computed. This can be seen in Figure~\ref{comparison}, where the two-grid convergence factors predicted by the analysis (denoted as $\rho$ and displayed as a solid line) are compared with those asymptotic convergence factors experimentally computed (represented by $\rho_h$ and depicted by using circles). To compute these latter, we consider a grid of size $256\times32$, and we use a $W-$cycle, a random initial guess and a zero right-hand side in order to avoid round-off errors.
We can see in the picture a very accurate prediction of the semi-algebraic mode analysis, making its use very useful for the analysis of the proposed multigrid waveform relaxation method.

Finally, we would like to show that the behavior of the proposed multigrid waveform relaxation is very satisfactory with respect to the number of time-levels considered. Since it is usually sufficient to analyze the behavior of the two-grid method to estimate the convergence of the multigrid method (see~\cite{TOS01}), in Table~\ref{lambdas_m}, we show the two-grid convergence factors provided by the analysis by considering a wide range of values of $M=2^k, \; k=5,\ldots,10$, together with the experimentally computed asymptotic convergence factors obtained by using the multilevel $W$-cycle with one smoothing step. As expected, the predicted two-grid convergence factors provide a very accurate estimation of the real asymptotic convergence of the method. These results are shown for different values of parameter $\lambda$, and a fixed value of the fractional order $\delta = 0.4$.

\begin{table}
\begin{center}
\begin{tabular}{ccccccc}
\hline
$\log_2 \lambda$ & $M=32$ & $M=64$ & $M=128$ & $M=256$ & $M=512$ & $M=1024$ \\
\hline
-8  & 0.004 & 0.005 & 0.006 & 0.008 & 0.010 & 0.013 \\
    & (0.004) & (0.005) & (0.006) & (0.009) & (0.012) & (0.014) \\
-6  & 0.016 & 0.018 & 0.023 & 0.028 & 0.036 & 0.045 \\
    & (0.017) & (0.018) & (0.027) & (0.033) & (0.041) & (0.051) \\
-4  & 0.054 & 0.061 & 0.072 & 0.085 & 0.098 & 0.110 \\
    & (0.055) & (0.065) & (0.079) & (0.088) & (0.091) & (0.112) \\
-2  & 0.116 & 0.120 & 0.125 & 0.128 & 0.132 & 0.134 \\
    & (0.122) & (0.125) & (0.135) & (0.137) & (0.137) & (0.137) \\
0   & 0.114 & 0.117 & 0.119 & 0.120 & 0.121 & 0.121 \\
    & (0.120) & (0.120) & (0.120) & (0.121) & (0.122) & (0.122) \\
2   & 0.054 & 0.058 & 0.060 & 0.062 & 0.062 & 0.063 \\
    & (0.057) & (0.063) & (0.063) & (0.063) & (0.063) & (0.063) \\
4   & 0.016 & 0.017 & 0.019 & 0.019 & 0.019 & 0.020 \\
    & (0.020) & (0.020) & (0.020) & (0.020) & (0.020) & (0.020) \\
6   & 0.004 & 0.004 & 0.005 & 0.005 & 0.005 & 0.005 \\
    & (0.005) & (0.005) & (0.005) & (0.005) & (0.005) & (0.005) \\
8   & 0.001 & 0.001 & 0.001 & 0.001 & 0.001 & 0.001 \\
    & (0.001) & (0.001) & (0.001) & (0.001) & (0.001) & (0.001) \\
\hline
\end{tabular}
\caption{Two-grid convergence factors predicted by the analysis together with the corresponding experimentally computed multilevel asymptotic convergence factors (between brackets) for different values of parameter $\lambda$ and for increasing number of time-steps, $M$, considering a fractional order $\delta = 0.4$. }
\label{lambdas_m}
\end{center}
\end{table}

\section{Fast implementation and computational cost}\label{sec:5}
In Algorithm~\ref{wrmg_1}, we observe that the most time-consuming part of the multigrid waveform
relaxation method is the calculation of the defect  and the smoothing step. The remaining components of the algorithm can be performed with a computational cost proportional to the number of unknowns. In the calculation of the residual, for each spatial grid-point a matrix-vector multiplication $T_M x$ is required for some vector $x$, where $T_M$ is the low-triangular matrix
\begin{equation}
T_M = \left(
\begin{array}{cccc}
 r_1 & 0 & \cdots & 0 \\
 r_2 & r_1  & \cdots & 0 \\
 \vdots & \ddots & \ddots & \vdots \\
 r_M & \cdots & r_2 & r_1
  \end{array}
\right),
\label{Toeplitz_matrix}
\end{equation}
with $r_i = \displaystyle \frac{\tau^{-\delta}}{\Gamma(2-\delta)} (d_i-d_{i-1}), i = 1, \ldots, M$, assuming $d_0 = 0$. \\
Moreover, the smoothing part involves the solution of triangular linear systems. The matrix ${\mathcal A}_{h,\tau}$ of the discrete system to solve  can be written as
${\mathcal A}_{h,\tau} = T_M \otimes I_N + A_h$, where $I_N$ denotes the identity matrix of order $N$, $A_h$ corresponds to the spatial discretization, $\otimes$ denotes the Kronecker product and  $T_M$ is the low-triangular matrix given
in \eqref{Toeplitz_matrix}. \\
In a standard implementation, the method would have a computational cost of at least  $O(NM^2)$ operations due to the matrix-vector multiplication $T_M x$ and the solution of the triangular systems in the smoothing part of the algorithm.
However, due to the special structure of matrix $T_M$, the proposed multigrid waveform relaxation method can be
implemented with a computational cost of $O(N M \log (M))$ operations with an storage cost for the system matrix of  $O(M)$. To see this, we discuss the following issues in the next subsections: a fast matrix-vector multiplication, a fast solution of the low-triangular systems, an efficient storage of matrix ${\mathcal A}_{h,\tau}$ and an estimation of the computational cost of the complete multigrid waveform relaxation method.

\subsection{An ${\mathbf {O(N M \log (M))}}$ calculation of the defect} To compute the residual in the multigrid waveform relaxation method, a matrix-vector multiplication ${\mathcal A}_{h,\tau} u$ is required. The matrix-vector multiplication corresponding to the spatial discretization can be calculated with a computational cost of $O(NM)$. Apart from this, for each spatial grid-point we have to perform a matrix-vector multiplication $T_M x$ for some vector $x$.  Notice that matrix $T_M$ is an $M\times M$ Toeplitz matrix, and therefore it can be embedded  into a $2M \times 2M$ circulant matrix
$C_{2M}$ in the following way:
\[
C_{2M} = \left(
\begin{array}{cc}
T_M & R_M \\
R_M & T_M
\end{array}
\right),
\]
where
$$
R_M = \left(
\begin{array}{ccccc}
 0 & r_{M} & r_{M-1} & \cdots & r_2 \\
 0 & 0  & r_M & \cdots & r_3 \\
 \vdots & \ddots & \ddots & \ddots & \vdots \\
 0 & \cdots & 0 & \hdots & r_{M} \\
 0 & 0 & 0  & \hdots & 0
  \end{array}
\right).
$$
Taking into account that
\[
 \left(
\begin{array}{cc}
T_M & R_M \\
R_M & T_M
\end{array}
\right)
 \left(
\begin{array}{c}
x \\
0
\end{array}
\right) =
 \left(
\begin{array}{c}
T_M x \\
*
\end{array}
\right),
\]
the matrix-vector multiplication is reduced to a circulant matrix-vector multiplication. It is known that a circulant matrix can be diagonalized by the Fourier matrix $F_{2M}$ as $C_{2M} = F_{2M}^* D_{2M} F_{2M}$,
where $D_{2M}$ is a diagonal matrix whose diagonal elements are the eigenvalues of $C_{2M}$. By taking the fast Fourier transform (FFT) of the first column of $C_{2M}$, we can determine  $D_{2M}$ in $O(M \log (M))$ operations. Once $D_{2M}$ is obtained, the multiplication $C_{2M} v$ for some vector $v$ can be calculated by using a couple of FFTs with $O(M \log (M))$ complexity. As this is the computational cost for each spatial grid-point, the product  ${\mathcal A}_{h,\tau} u$ can be performed in $O(N M \log (M))$ operations.

\subsection{An ${\mathbf {O(N M \log (M))}}$ implementation of the smoothing procedure}
Other of the most consuming components of the multigrid waveform relaxation method for solving the time-fractional diffusion equation is the relaxation step, since dense low triangular systems must be solved. In the particular case of discrete problem \eqref{discrete_model_IVP_1}, for each spatial grid-point  we need to solve a system of $M$ equations of the type $(T_M + 2/h^2 I_M) x = b$ for some known vector $b$.  Due to the Toeplitz-structure of the matrix, the solution of the system can be obtained  in
$O(M \log (M))$ operations by using well-developed algorithms for the inversion of triangular Toeplitz matrices. With the inverse matrix obtained, which is again a Toeplitz matrix, the solution of the system is obtained by a matrix-vector multiplication with complexity of $O(M \log (M))$ operations by using the algorithm described in the previous subsection. Classical algorithms for the inversion of triangular Toeplitz matrices with complexity $O(M \log (M))$ include the Bini's algorithm \cite{bini}, its revised version \cite{lin_ching}, and the divide and conquer method \cite{commenges,morf}. In our implementation, we have chosen the latter,  which is briefly described  to make this work more self-contained. A low-triangular Toeplitz matrix $T_M$, with $M=2^p, p>1$,  can be partitioned as follows
\[
T_{M} = \left(
\begin{array}{cc}
T_{M/2} &  \\
P_{M/2} & T_{M/2}
\end{array}
\right),
\]
where $T_{M/2}$ and $P_{M/2}$ are Toeplitz matrices of order $2^{p-1}$. Based on this partition, it is easy
to see that the inverse of matrix $T_{M}$ can be written as
\[
T_M^{-1} = \left(
\begin{array}{cc}
T_{M/2}^{-1} &  \\
-T_{M/2}^{-1} P_{M/2} T_{M/2}^{-1}& T_{M/2}^{-1}
\end{array}
\right).
\]
This expression gives us a recurrent method to calculate the inverse of matrix $T_M$.
Since the inverse of this matrix is Toeplitz, it is enough to calculate its first column.
Given a small number $p_0$, we compute the inverse of the submatrix $T_{2^{p_0}}$ by the
forward substitution method, for instance. Then we subsequently apply the recurrent formula to
compute the inverse of $T_{M}$ in $p - p_0$ steps. On each step 
the first column of the Toeplitz matrix $-T_{M/2}^{-1} P_{M/2} T_{M/2}^{-1}$ is required, which can be calculated by FFTs. The total computational cost of the  smoothing algorithm is therefore only $O(N M \log (M))$
at each iteration step. Moreover, since we need to solve several triangular systems
with the same matrix but different right-hand sides, the first column of the inverse matrix can be computed
a priori.
\subsection{Storage cost and computational complexity of the multigrid waveform relaxation method}
The non-local nature of the fractional derivatives results in a dense coefficient matrix yielding a bottleneck
for the traditional numerical methods for fractional diffusion problems which require $O(M^2)$ units of
storage. Due to the Toeplitz-structure of matrix $T_M$ the memory requirement for the storage of the coefficient matrix can be significantly reduced
to $O(M)$, since to perform all the calculations in our algorithm, we only need to store its first column. \\
We consider a grid-hierarchy $G_0, G_1, \ldots, G_l$, where $G_k := G_{h_k,\tau}$ and $h_0>h_1>\ldots>h_l$.
It is well-known that the computational work ${\cal W}_l$ per $V-$cycle on a grid $G_l$ is given by \cite{Stu_Tro}
\[
{\cal W}_l = \sum_{k=1}^l {\cal W}_k^{k-1} + {\cal W}_0,
\]
where ${\cal W}_k^{k-1}$ is the computational work of a two-grid cycle excluding
the work needed to solve the defect equation on $G_k$, and ${\cal W}_0$ is the work needed to
compute the exact solution on the coarsest grid $G_0$. In the computational work ${\cal W}_k^{k-1}$,
it is included the cost of a smoothing iteration, the calculation of the defect and its transfer to $G_{k-1}$,
and the interpolation of the correction to $G_k$ and  its addition to the previous approximation. From the
previous subsections,  we can estimate that the computational cost of a two-grid cycle is
${\cal W}_k^{k-1} = O(N_k M \log (M))$ and on the coarsest grid  ${\cal W}_0 = O(M \log (M))$, where $N_k$ is the number of spatial grid-points on the grid $G_k$ and $M$ is the number of time steps. Therefore, we can say that the  computational cost of a $V-$cycle on level $l$ is roughly
\[
{\cal W}_l = (1+ \frac{1}{2} + \frac{1}{2^2} + \ldots + \frac{1}{2^l}) O(N_l M \log (M)) = O(N_l M \log (M)),
\]
Thus, since the $V-$ cycle converges in a small number of iterations independent of the number of unknowns,
the total computational cost for solving the time-fractional problem by the multigrid waveform relaxation method
is roughly $O(N_l M \log (M))$.

\section{Extension to 2D}\label{sec:6}

\setcounter{section}{6}

This section is devoted to the extension of the presented methodology to problems with two spatial dimensions. \\

\noindent {\bf Model problem and discretization.}
We consider the two-dimensional time-frac\-tional diffusion equation as model problem, that is,
\begin{eqnarray}
D_t^{\delta} u - \Delta u &=& f(x,y,t),\quad (x,y)\in\Omega \subset {\mathbb R}^2,\; t>0, \label{2D_model_BVP_1}\\
u(x,y,t) &=& 0,\quad (x,y)\in \partial\Omega, \; t>0,\label{2D_model_BVP_2}\\
u(x,y,0) &=& g(x,y), \quad (x,y)\in\overline{\Omega},\label{2D_model_BVP_3}
\end{eqnarray}
where $\Delta$ denotes the two-dimensional Laplace operator, $\Omega$ is a square domain of length $L$, $\partial \Omega$ is its boundary and $\overline{\Omega} = \Omega \cup \partial \Omega$. $D_t^{\delta}$ denotes again the Caputo fractional derivative,
\begin{equation}\label{caputo_bis}
D_t^{\delta} u (x,y,t):= \left[\frac{1}{\Gamma(1-\delta)}\int_{0}^{t} (t-s)^{-\delta}\frac{\partial u(x,y,s)}{\partial s} ds\right],\quad (x,y)\in\Omega,\; t\geq 0.
\end{equation}
Let us consider a uniform grid $G_{h,\tau} = G_h\times G_{\tau}$, with
\begin{equation}\label{spatial_mesh_2d}
G_h = \left\{(x_n,y_l)\,| \, x_n = nh,\, y_l = lh, \, n,l=0,1,\ldots,N+1\right\},
\end{equation}
where $h=\displaystyle \frac{L}{N+1}$, and with $G_{\tau}$ given as in~\eqref{temporal_mesh}. The nodal approximation to the solution at each grid point $(x_n,y_l,t_m)\in G_{h,\tau}$ is denoted by $u_{n,l,m}$.\\
Standard central finite differences are used again to approximate the spatial derivatives, whereas the Caputo fractional derivative is discretized as
\begin{equation}
D_M^{\delta}u_{n,l,m}:=\frac{\tau^{-\delta}}{\Gamma(2-\delta)}\left[d_1 u_{n,l,m} - d_m u_{n,l,0} + \sum_{k=1}^{m-1}(d_{k+1}-d_k)u_{n,l,m-k}\right],\label{L1_scheme_2d}
\end{equation}
where coefficients $d_k$ are identically defined as in Section~\ref{sec:2}.\\
This results in the following discrete problem
\begin{eqnarray}
D_M^{\delta} u_{n,l,m}- \Delta_h u_{n,l,m} &=& f(x_n,y_l,t_m), 1\leq n,l\leq N,\; 1\leq m\leq M,\label{discrete_model_IVP_2d_1}\\
u_{n,l,m} &=& 0,\;  (x_n,y_l)\in\partial \Omega\cap G_{h,\tau},\; 0<m\leq M,\label{discrete_model_IVP_2d_2}\\
u_{n,l,0} &=& g(x_n,y_l), 0\leq n,l\leq N+1,\label{discrete_model_IVP_2d_3}
\end{eqnarray}
where
$$\Delta_h u_{n,l,m} = \frac{u_{n+1,l,m}+u_{n,l+1,m}-4u_{n,l,m}+u_{n-1,l,m}+u_{n,l-1,m}}{h^2}.$$

\noindent {\bf Multigrid waveform relaxation in 2D.}
Regarding the solver for the considered two-dimensional time-fractional model problem~\eqref{2D_model_BVP_1}-\eqref{2D_model_BVP_3}, a red-black Gauss-Seidel waveform relaxation can be defined, after discretizing in space, as follows,
\begin{eqnarray}\label{red-black_2d}
&&D_t^{\delta}u_{n,l}^k(t) +\frac{4}{h^2}u_{n,l}^k(t) = \frac{1}{h^2}\left(u_{n-1,l}^{k-1}(t)+u_{n,l-1}^{k-1}(t)+u_{n+1,l}^{k-1}(t)+u_{n,l+1}^{k-1}(t) \right)\nonumber\\
&& \qquad \qquad \qquad \qquad \qquad \qquad+f_{n,l}(t),\; \hbox{if} \; n+l\; \hbox{is even},\\
&&D_t^{\delta}u_{n,l}^k(t)-\frac{1}{h^2}\left(u_{n-1,l}^{k}(t)+u_{n,l-1}^{k}(t)-4u_{n,l}^k(t)+u_{n+1,l}^{k}(t)+u_{n,l+1}^{k}(t)\right)\nonumber \\
&& \qquad \qquad \qquad \qquad \qquad \qquad= f_{n,l}(t),\; \hbox{if} \; n+l\; \hbox{is odd}.
\end{eqnarray}
Thus, the fully discrete problem given in~\eqref{discrete_model_IVP_2d_1}-\eqref{discrete_model_IVP_2d_3} can be solved by using an extension of the multigrid waveform relaxation algorithm proposed in Section~\ref{sec:3}. In this case, the method combines a two-dimensional coarsening strategy in the space variables and again a line-in-time smoother based on the red-black Gauss-Seidel waveform relaxation, that is, the lines in time are visited following a red-black or chessboard manner. Regarding the inter-grid transfer operators, the standard two-dimensional full-weighting restriction and bilinear interpolation are considered.\\

\noindent {\bf Semi-algebraic mode analysis in 2D.}
The semi-algebraic mode analysis presented in Section~\ref{sec:4} can be also extended to study the convergence of the proposed multigrid waveform relaxation method. For this analysis, very little has to be changed from the theory developed in Section~\ref{sec:4}. The infinite grid ${\mathcal G}_h$ is defined as the extension of the spatial mesh given in~\eqref{spatial_mesh_2d}, and then, the grid-functions defined on such a grid can again be expressed as formal linear combinations of the Fourier components which in this case are given by the product of two complex exponential functions, i.e. $\varphi_h({\boldsymbol \theta},{\mathbf x}) = e^{\imath {\boldsymbol \theta} \cdot {\mathbf x}} = e^{\imath \theta_x \, x}e^{\imath \theta_y\, y},$ where ${\boldsymbol \theta} = (\theta_x,\theta_y)\in{\boldsymbol \Theta}_h = (-\pi/h,\pi/h]\times (-\pi/h,\pi/h]$, and which form the new Fourier space. In the two-dimensional spatial case, it is well-known that the Fourier space is decomposed in four-dimensional subspaces
$${\mathcal F}^4 ({\boldsymbol \theta}) = \hbox{span}\left\{\varphi_h({\boldsymbol \theta}^{00},\cdot), \, \varphi_h({\boldsymbol \theta}^{11},\cdot), \, \varphi_h({\boldsymbol \theta}^{10},\cdot), \, \varphi_h({\boldsymbol \theta}^{01},\cdot) \right\},$$
generated by four Fourier modes associated with one low frequency ${\boldsymbol \theta} = {\boldsymbol \theta}^{00}\in {\Theta}_{2h} = [-\pi/2h,\pi/2h)^2$ and three high frequencies ${\boldsymbol \theta}^{11}$, ${\boldsymbol \theta}^{10}$ and ${\boldsymbol \theta}^{01}$ such that,
$${\boldsymbol \theta}^{\alpha_1,\alpha_2} = {\boldsymbol \theta}^{00} - (\alpha_1\, \hbox{sign}(\theta_1^{00})\pi, \, \alpha_2\,\hbox{sign}(\theta_2^{00})\pi), \; \alpha_1,\alpha_2\in\{0,1\},$$
which are coupled on the coarse grid by the aliasing effect. \\
Similarly as in the one-dimensional case, the semi-algebraic mode analysis in 2D is based on a two-dimensional spatial local Fourier analysis combined with an exact analysis in time. In this way, the resulting Fourier representations of the smoothing, coarse-grid and two-grid operators are $4M\times 4M$ matrices.\\

\begin{figure}[htb]
\begin{center}
\includegraphics[width = 0.8\textwidth]{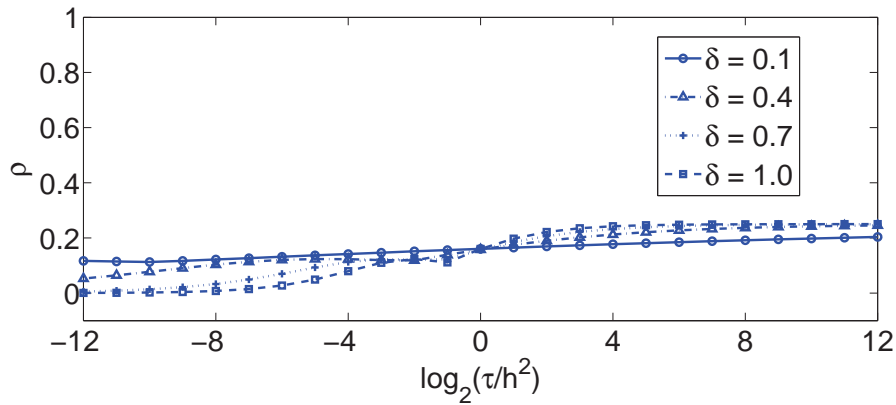}
\caption{Two-grid convergence factors predicted by the analysis for different values of parameter $\lambda = \tau/h^2$ and various fractional orders $\delta$.}
\label{delta_lambda_star_2d}
\end{center}
\end{figure}

\begin{figure}[htb]
\begin{center}
\includegraphics[width = 0.8\textwidth]{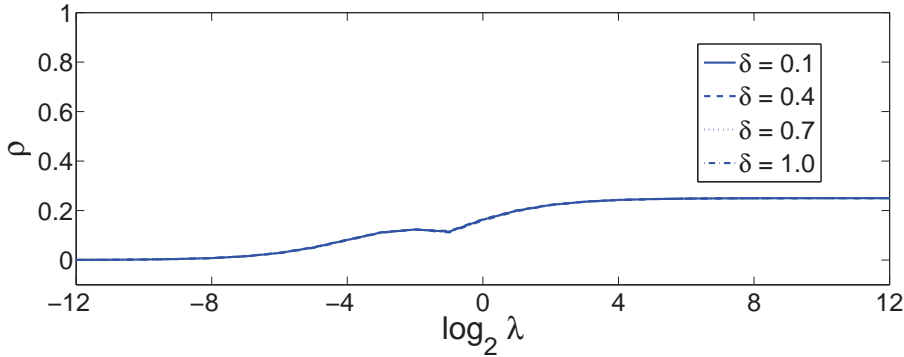}
\caption{Two-grid convergence factors predicted by the analysis for different values of parameter $\lambda = \tau^{\delta} \Gamma(2-\delta)/h^2$ and various fractional orders $\delta$.}
\label{delta_lambda_2d}
\end{center}
\end{figure}

\begin{figure}[htb]
\begin{center}
\includegraphics[width = 0.8\textwidth]{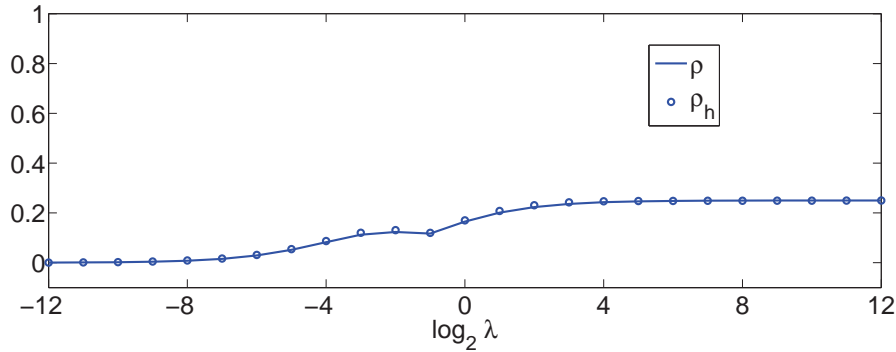}
\caption{Comparison between the two-grid convergence factors predicted by the analysis ($\rho$) and the asymptotic convergence factor of a $W-$cycle experimentally computed ($\rho_h$), for different values of parameter $\lambda = \tau^{\delta} \Gamma(2-\delta)/h^2$ and fractional order $\delta = 0.4$.}
\label{comparison_2d}
\end{center}
\end{figure}

\noindent {\bf Analysis results.}
Next, we present some results obtained by using the semi-algebraic analysis.
Similarly as we saw for the 1D model problem, if we analyze the convergence of the method depending on parameter $\tau/h^2$, although the convergence rates are bounded by $0.25$ for all cases, we do not obtain a $\delta-$independent convergence for a fixed value of $\tau/h^2$. This can be seen in Figure~\ref{delta_lambda_star_2d}.
However, we can show that the obtained multigrid convergence becomes robust for any value of $\delta$ with respect to parameter $\lambda = \tau^{\delta} \Gamma(2-\delta)/h^2$. This is shown in Figure~\ref{delta_lambda_2d}, where $M = 32$ time-levels have been considered, and the two-grid convergence factors predicted by the analysis for one smoothing step are shown for different values of parameter $\lambda$ and for different fractional orders $\delta$. Notice that the graphs corresponding to the different values of $\delta$ are almost indistinguishable, and for any value of $\lambda$ the multigrid convergence results very satisfactory.
These results can be confirmed with the asymptotic convergence factors experimentally computed. In particular, for $\delta = 0.4$, we show this comparative in Figure~\ref{comparison_2d}, where the two-grid convergence factors predicted by the semi-algebraic mode analysis are displayed together with the asymptotic convergence rates computed by using a $W(1,0)-$multigrid waveform relaxation algorithm on a fine grid of size $256\times256\times 32$. Again, a random initial guess and a zero right-hand side are used to perform these calculations.
Similar pictures can be obtained for other fractional orders $\delta$.
We can observe a very accurate match between the analysis results and the rates experimentally obtained.

\noindent \textbf{Remark.} In Figures~\ref{comparison} (1d case) and ~\ref{comparison_2d} (2d case), it is observed a different behavior of the multigrid method when $\lambda$ becomes big enough, that is, in the limit case of the steady problem. In that case, it is well-known that the multigrid method based on a red-black smoother is an exact solver in the one-dimensional case whereas for a two-dimensional diffusion problem the convergence rate is about $0.25$ for a $W-$cycle with one smoothing step (see~\cite{TOS01}).

\section{Numerical results}\label{sec:7}

In this section, we consider three different numerical experiments to illustrate the efficiency of the proposed multigrid waveform relaxation method for solving the time-fractional heat equation. For all cases we perform $V-$cycles since they provide similar convergence rates to $W-$cycles and therefore a more efficient multigrid method is obtained. We will start solving both one- and two-dimensional linear problems and finally we will solve a non-linear one-dimensional problem. All numerical computations were carried out using MATLAB.\\

\noindent{\bf One-dimensional linear time-fractional heat equation.}
We show the efficient performance of the proposed multigrid waveform relaxation for a problem which considers reasonably general and realistic hypotheses on the behavior of the solution near the initial time. In particular, we consider a problem whose solution is smooth away from the initial time ($t=0$) but it has a certain singular behavior at $t=0$ presenting a boundary layer. The theoretical convergence analysis of the considered finite difference discretization has been deeply studied in~\cite{martin1}. Here, we will show that the convergence of the WRMG is satisfactory for this representative model problem.

\begin{figure}[htb]
\begin{center}
\includegraphics[width = 0.8\textwidth]{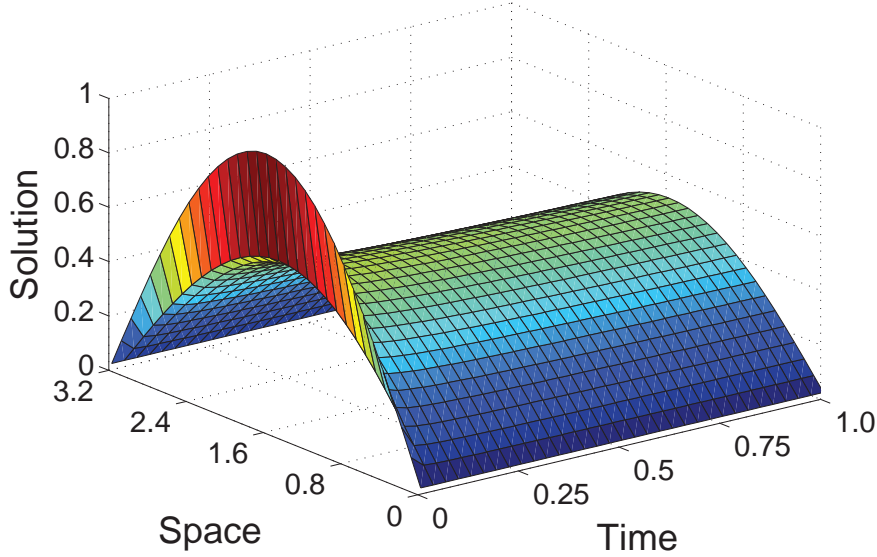}
\caption{Analytical solution $u(x,t)$ of the first test problem, for fractional order $\delta=0.1$.}
\label{solution_picture}
\end{center}
\end{figure}

\begin{figure}[htb]
\begin{center}
\includegraphics[width = 0.8\textwidth]{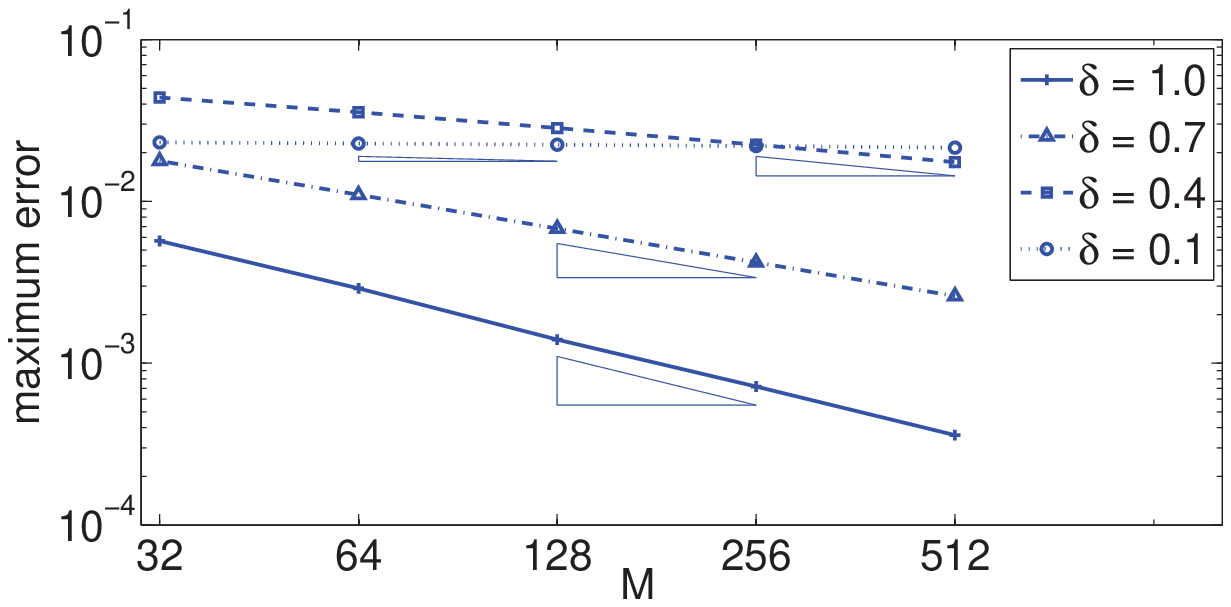}
\caption{Reduction of the maximum errors obtained for four different values of $\delta$, for the first test problem.}
\label{error_reduction}
\end{center}
\end{figure}

We consider problem~\eqref{model_IVP_1}-\eqref{model_IVP_3} defined on a domain $[0,\pi]\times [0,1]$, with a zero right-hand side ($f(x,t) = 0$) and an initial condition $g(x) = \sin\, x$. Then, function $u(x,t) = E_{\delta}(-t^{\delta})\sin\, x$, where $E_{\delta}:{\mathbb R} \rightarrow {\mathbb R}$ is given by
$$E_{\delta}(z) := \sum_{k=0}^{\infty}\frac{z^k}{\Gamma(\delta \, k + 1)},$$
satisfies our initial-boundary value problem~\cite{Luchko, martin1}. In Figure~\ref{solution_picture}, we can observe the sharpness of the analytical solution near the initial time, where a boundary layer appears.
In~\cite{martin1}, it is proved rigorously that for ``typical'' solutions of~\eqref{model_IVP_1}-\eqref{model_IVP_3} (no excessive smooth solutions) a rate of convergence of ${\mathcal O}(h^2+\tau^{\delta})$ is obtained.  This is shown in Figure~\ref{error_reduction} for four different values of $\delta$, where the maximum errors between the analytical and the numerical solution are displayed for various numbers of time-steps $M$ and assuming a sufficiently fine spatial grid.
It can be seen that the slopes of the obtained graphs match with the expected convergence rates.
For small values of $\delta$, a very fine temporal mesh would be required to attain the asymptotic rate of convergence, and this is the case of $\delta = 0.1$ in the picture where a slow convergence of the rates to the expected asymptotic rate of $0.1$ is observed.

Next, we show the independence of the convergence of the proposed multigrid waveform relaxation method with respect to the discretization parameters. For this purpose, we consider different values of the fractional order $\delta$ and different grid-sizes varying from $128\times 128$ to $2048\times 2048$ doubling the mesh-size in both spatial and temporal dimensions.
In Table~\ref{table_it_1d} we display the number of WRMG iterations necessary to reduce the initial residual in a factor of $10^{-10}$, together with the mean convergence factors and the corresponding CPU time, when considering a $V(0,1)-$cycle. We can observe that the performance of the $V-$cycle is also satisfactory for any value of $\delta$ and for increasing mesh-sizes, as it was already seen for the $W-$cycle in the analysis results section. Moreover, we choose only one post-smoothing step since this approach provides much better convergence factors than a $V(1,0)-$cycle. Taking into account these considerations, we observe from Table~\ref{table_it_1d} a robust convergence of the considered WRMG.\\
\begin{table}
\begin{center}
\begin{tabular}{cccccc}
\hline
$\delta$ & $128\times 128$ & $256\times 256$ & $512\times 512$ & $1024\times 1024$ & $2048\times 2048$ \\
\hline
0.1 & 8 (0.03) 0.54s & 8 (0.03) 1s & 8 (0.03) 2.96s & 8 (0.03) 10.54s & 7 (0.03) 36.16s \\
0.4 & 7 (0.03) 0.49s & 7 (0.03) 0.91s & 7 (0.03) 2.60s & 7 (0.03) 9.31s & 7 (0.03) 36.34s \\
0.7 & 7 (0.04) 0.47s & 7 (0.04) 0.90s & 7 (0.04) 2.54s & 7 (0.04) 9.15s & 7 (0.04) 36.10s \\
1.0 & 7 (0.05) 0.46s & 7 (0.05) 0.88s & 7 (0.05) 2.54s & 6 (0.05) 7.96s & 6 (0.05) 30.69s \\
\hline
\end{tabular}
\caption{Number of $V(0,1)-$WRMG iterations necessary to reduce the initial residual in a factor of $10^{-10}$ for different fractional orders $\delta$ and for different grid-sizes. The corresponding average convergence factors (between brackets) and the CPU times in seconds are also included.}
\label{table_it_1d}
\end{center}
\end{table}

\noindent{\bf Two-dimensional linear time-fractional heat equation.}
The second numerical experiment deals with the solution of a two-dimensional linear time-fractional diffusion problem. We consider the following model problem defined on the spatial domain $\Omega=(0,2)\times (0,2)$
\begin{eqnarray}
D_t^{\delta} u - \Delta u &=& f(x,y,t),\quad (x,y)\in\Omega, \; t>0, \label{2D_model_BVP_1_ex}\\
u(x,y,t) &=& 0,\quad (x,y)\in \partial\Omega, \; t>0,\label{2D_model_BVP_2_ex}\\
u(x,y,0) &=& 0, \quad (x,y)\in\overline{\Omega},\label{2D_model_BVP_3_ex}
\end{eqnarray}
where $$f(x,y,t) = \left(\frac{2t^{2-\delta}}{\Gamma(3-\delta)}+\left(1+\frac{\pi^2}{2}\right)t^2\right)\sin \frac{\pi\,x}{2}\sin\frac{\pi\, y}{2},$$
in the way that the analytic solution of the problem is $$u(x,y,t) = t^2 \sin\frac{\pi\, x}{2}\sin \frac{\pi\,y}{2}.$$
We consider the multigrid waveform relaxation method described in Section~\ref{sec:6}, by using a $V(1,1)-$cycle. This choice is based on the semi-algebraic mode analysis results presented in Section~\ref{sec:6}. Due to the difference of the behavior of the method between the one- and two-dimensional problems, we have chosen two-smoothing steps to perform the calculations in this test case. \\
In Table~\ref{table_it_2d} we display the number of WRMG iterations necessary to reduce the initial residual in a factor of $10^{-10}$ for different grid-sizes varying from $32\times 32\times 32$ to $256\times 256\times 256$ and for different values of the fractional order $\delta$. We can observe that the convergence of the proposed multigrid waveform relaxation is very robust with respect to the considered parameters. In the table, we also show the mean convergence factors and the corresponding CPU times. 
\begin{table}
\begin{center}
\begin{tabular}{ccccc}
\hline
$\delta$ & $32\times 32\times 32$ & $64\times 64\times 64$ & $128\times 128\times 128$ & $256\times 256\times 256$  \\
\hline
0.1 & 12 (0.10) 2.46s & 12 (0.10) 10.31s & 12 (0.11) 55.52s & 12 (0.11) 349.98s \\
0.4 & 12 (0.09) 2.51s & 12 (0.10) 10.45s & 12 (0.11) 55.86s & 12 (0.11) 348.29s \\
0.7 & 11 (0.09) 2.27s & 12 (0.10) 10.31s & 12 (0.11) 55.63s & 12 (0.11) 344.57s \\
1.0 & 11 (0.09) 2.29s & 11 (0.10) 9.68s & 12 (0.11) 55.73s & 12 (0.11) 346.44s \\
\hline
\end{tabular}
\caption{Number of $V(1,1)-$WRMG iterations necessary to reduce the initial residual in a factor of $10^{-10}$, together with the corresponding average convergence factors (between brackets) and the CPU times in seconds, for different fractional orders $\delta$ and for different grid-sizes.}
\label{table_it_2d}
\end{center}
\end{table}
We can observe a very satisfactory convergence in all cases, making the multigrid waveform relaxation method a good choice for an efficient solution of the time-fractional two-dimensional heat equation. \\

\noindent{\bf One-dimensional nonlinear problem.}
The last numerical experiment is devoted to deal with a nonlinear problem which appears in the modeling of anomalous diffusion in porous media~\cite{Hilfer, SIAM_plociniczak}. We consider the following time-fractional partial differential equation
\begin{equation}\label{anomalous_diffusion}
D_t^{\delta} u = \frac{\partial}{\partial x}\left(D(u)\frac{\partial u}{\partial x}\right)+c\frac{\partial u}{\partial x} + f(x,t),
\end{equation}
where $D_t^{\delta}$ denotes again the Caputo fractional derivative operator with $0<\delta<1$, and $f(x,t)$ represents a source term. In this test problem we assume homogeneous Dirichlet boundary conditions and a zero initial condition. \\
Choosing $c=0$, model problem~\eqref{anomalous_diffusion} has been used to describe the moisture distribution in construction materials~\cite{paper_plociniczak} for example, whereas if the convective term is included it is used to describe transport models for single-phase gas through tight rocks~\cite{shale_gas} or in groundwater hydrology~\cite{Baeumer}.
For the discretization of problem~\eqref{anomalous_diffusion}, we consider again a uniform grid in space and time with step-sizes $h$ and $\tau$, respectively. The fractional temporal derivative is discretized as previously by using the L1 scheme (see~\eqref{L1_scheme2}). Regarding the spatial discretization, in an interior point $(x_n,t_m)$ the diffusion term is approximated by
\begin{equation}\label{spatial_discretization}
\frac{1}{h}\left[a_{n+1/2,m}\frac{u_{n+1,m}-u_{n,m}}{h} - a_{n-1/2,m}\frac{u_{n,m}-u_{n-1,m}}{h}\right],
\end{equation}
where $a_{n\pm 1/2,m} = \displaystyle\frac{1}{2}\left[D(u_{n\pm1,m})+D(u_{n,m})\right]$, and for the convective term a standard upwind scheme is considered.\\
For the solution of the resulting discrete problem, we propose a nonlinear multigrid waveform relaxation method, that is, the well-known waveform relaxation FAS method. This algorithm is easily derived from the standard FAS method~\cite{TOS01} for solving elliptic equations. For a detailed description of the proposed algorithm we refer the reader to the book~\cite{Vandewalle_book}. A nonlinear Gauss-Seidel waveform relaxation with a red-black ordering is considered, together with standard transfer-grid operators. Again, a $V(0,1)-$cycle is chosen to perform the calculations. \\
In Table~\ref{table_nonlinear}, we show the convergence of the proposed algorithm for the case of $D(u) = 1+u^2$, $c=1$ and $f(x,t) = 1$, and for different values of the fractional order $\delta$. In particular, we display the number of iterations required to reduce the maximum initial residual by a factor of $10^{-10}$ for different grid-sizes and the corresponding mean convergence factors (between brackets).
\begin{table}
\begin{center}
\begin{tabular}{ccccc}
\hline
$\delta$ & $32\times 32$ & $64\times 64$ & $128\times 128$ & $256\times 256$ \\
\hline
0.1 & 11 (0.09) & 11 (0.10) & 11 (0.10) & 12 (0.10) \\
0.4 & 11 (0.09) & 11 (0.10) & 12 (0.10) & 12 (0.10) \\
0.7 & 11 (0.09) & 11 (0.10) & 12 (0.10) & 12 (0.10) \\
1.0 & 11 (0.10) & 12 (0.10) & 12 (0.10) & 12 (0.10) \\
\hline
\end{tabular}
\caption{Number of $V(0,1)-$iterations of the waveform relaxation FAS method required to reduce the initial residual in a factor of $10^{-10}$ for different fractional orders $\delta$ and for different grid-sizes, together with the corresponding mean convergence factors (between brackets).}
\label{table_nonlinear}
\end{center}
\end{table}
From the results in Table~\ref{table_nonlinear}, we can conclude that the waveform relaxation FAS method shows a similar behavior as the linear multigrid waveform relaxation method for the time-fractional diffusion problems.

\section{Conclusions}\label{sec:8}

A multigrid waveform relaxation method has been proposed for solving the time-fractional heat equation. The convergence of this method has been studied by a suitable semi-algebraic mode analysis, which combines a classical exponential Fourier analysis in space with an algebraic computation in time. The results of this analysis show the efficiency and robustness of the proposed algorithm for the solution of the considered problem for different fractional orders.
The proposed method has a computational cost of $O(N M \log(M))$ operations, where $M$ is the number of time steps and $N$ is the number of spatial grid points.
Moreover, three numerical experiments confirm the good behavior of the WRMG method. In particular a linear one-dimensional representative problem, a linear two-dimensional model problem and a nonlinear one-dimensional problem with applications in porous media are efficiently solved in this work.

%
\section*{Acknowledgments}
The authors thank the referees for their valuable comments and suggestions which helped to improve the paper.

\bibliographystyle{plain}
\bibliography{references}
\end{document}